\documentclass[12pt]{amsart}
\usepackage{amsmath,amssymb,amsthm}
\usepackage[]{latexsym}
\newtheorem{defn0}{Definition}[section]
\newtheorem{prop0}[defn0]{Proposition}
\newtheorem{thm0}[defn0]{Theorem}
\newtheorem{lemma0}[defn0]{Lemma}
\newtheorem{corollary0}[defn0]{Corollary}
\newtheorem{example0}[defn0]{Example}
\newtheorem{remark0}[defn0]{Remark}
\newtheorem{conjecture0}[defn0]{Conjecture}
\newtheorem{notation0}[defn0]{Notation}

\newenvironment{definition}{\begin{defn0}\rm}{\end{defn0}}
\newenvironment{proposition}{\begin{prop0}}{\end{prop0}}
\newenvironment{theorem}{\begin{thm0}}{\end{thm0}}

\newenvironment{corollary}{\begin{corollary0}}{\end{corollary0}}

\newenvironment{remark}{\begin{remark0}\rm}{\end{remark0}}

\newcommand{\cO}{{\mathcal O}}

\newcommand{\Gal}{{\mathrm {Gal }}}
\newcommand{\disc}{{\mathrm {disc }}}

\newcommand{\M}{\mathrm{M }}

\newcommand{\GL}{{\mathrm \mathbb{GL}}}
\newcommand{\End}{{\mathrm{End }}}

\newcommand{\Fr}{{\mathrm {Frob }}}

\newcommand{\Z}{{\mathbb Z}}
\newcommand{\Q}{{\mathbb Q}}

\newcommand{\F}{{\mathbb F}}

\newcommand{\TT}{{\mathbb T}}

\newcommand{\qbar}{{\overline {\Q }}}

\include{thebibliography}

\begin{document}

\title{On abelian surfaces with potential quaternionic multiplication}

\author{Luis V. Dieulefait, Victor Rotger}

\footnote{Supported by a postdoctoral fellowship from the European
Research Network "Arithmetic Algebraic Geometry"; Partially
supported by DGICYT Grant BFM2003-06768-C02-02.}

\address{Centre de Recerca Matematica, Apartat 50, E-08193, Bellaterra,
Spain; Universitat Polit\`{e}cnica de Catalunya, Departament de
Matem\`{a}tica Aplicada IV (EUPVG), Av. Victor Balaguer s/n, 08800
Vilanova i la Geltr\'{u}, Spain.}

\email{Ldieulefait@crm.es, vrotger@mat.upc.es}



\begin{abstract}
An abelian surface $A$ over a field $K$ has potential quaternionic
multiplication if the ring $\End _{\bar K}(A)$ of geometric
endomorphisms of $A$ is an order in an indefinite rational
division quaternion algebra. In this brief note, we study the
possible structures of the ring of endomorphisms of these surfaces
and we provide explicit examples of Jacobians of curves of genus
two which show that our result is sharp.
\end{abstract}

\maketitle

\section{The ring of endomorphisms of an abelian surface with
potential quaternionic multiplication}

\begin{definition}\label{def}

Let $K$ be a field and let $\bar K$ be a separable closure of $K$.
An abelian surface $A$ over $K$ has potential quaternionic
multiplication if $\End _{\bar K}(A)$ is an order $\cO $ in an
indefinite quaternion division algebra $B$ over $\Q $.

\end{definition}

In the literature, it is often required that abelian surfaces with
quaternionic multiplication over a field $K$ have all the
endomorphisms defined over the base field. Under these hypothesis,
these abelian surfaces are also sometimes called fake elliptic
curves, since their arithmetic bear a strong analogy with the
arithmetic of elliptic curves. We refer the reader to \cite{Bu},
\cite{Jo1}, \cite{Jo2}, \cite{DiRo}, \cite{Oh} and \cite{Ro}.
However, in our definition we allow the endomorphisms to be
defined over a field extension $L/K$, since this is precisely the
phenomenon we wish to study and natural examples of them arise in
the modular setting.

When working with Shimura curves as moduli spaces of abelian
surfaces with quaternionic multiplication, one often considers
abelian surfaces that {\em contain} a quaternion order. Hence, in
our definition we are missing products of elliptic curves with
complex multiplication. Indeed, according to Definition \ref{def},
our abelian surfaces are absolutely simple. This is not relevant
for our purposes, since the arithmetic of elliptic curves with
complex multiplication is better known.

Let us recall that an order $\cO $ in a quaternion algebra is {\em
hereditary} if all its one-sided ideals are projective. An order
$\cO $ is hereditary if and only if the reduced discriminant
$D=\disc (\cO )$ is a square-free integer.

In \cite{DiRo}, the authors studied the possible fields of
definition of the quaternionic multiplication on an abelian
surface and the possible structures for the integral quaternion
order $\End _{\bar K}(A)$. In the following statement we recall
and strengthen these results in the case that the base field is
$\Q $.

\begin{theorem}\label{ModOrder}

Let $A/\Q $ be an abelian surface with potential quaternionic
multiplication by an hereditary order $\cO $ of discriminant $D$
in a quaternion algebra $B$. Let $L/\Q $ be the minimal field of
definition of the quaternionic multiplication on $A$ and $F=\End
_{\Q }(A) \otimes \Q$. Then, either

\begin{itemize}
\item $L$ is an imaginary quadratic field and $F=\Q (\sqrt {m})$
is a real quadratic field such that $B \simeq (\frac {-D \delta ,
m}{\Q })$ for any possible degree $\delta $ of a polarization on
$A/\Q $, or

\item $L$ is a purely imaginary dihedral extension of $\Q $ of
degree $[L:\Q ]\geq 4$ and $F = \Q $. If $A$ is the Jacobian
variety of a curve $C/\Q $ of genus two, then $[L:\Q ]= 4$.
\end{itemize}

\end{theorem}

{\em Proof. } As it is shown in \cite{DiRo}, the field extension
$L/\Q $ is Galois and $\Gal (L/\Q )\simeq C_n$ or $D_n$, $n\leq
6$, where $C_n$ is the cyclic group of order $n$ and $D_n$ denotes
the dihedral group of $2 n$ elements. Moreover, it is shown that
if $\End _{\Q }(A)\otimes \Q $ is a totally real number field over
$\Q $, then $L/\Q $ is necessarily dihedral. In addition, if $A$
admits a polarization over $\Q $ of degree not equal to $D$ nor $3
D$ up to squares, then $|\Gal (L/\Q )|=4$. This is for instance
the case when $A = \mbox{Jac }(C)$ where $C/\Q $ is a curve of
genus two. The statement above now follows from the following
proposition.

\begin{proposition}

Let $A/\Q $ be an abelian surface over $\Q $ which is simple over
$\qbar $. Then $\End _{\Q }(A)$ is either $\Z $ or an order in a
real quadratic field.

\end{proposition}

{\em Proof.} By Albert's classification of involuting division
algebras and Shimura's work on endomorphisms of abelian varieties,
$\End _{\qbar }(A)\otimes \Q $ is either $\Q $, a real quadratic
field or a totally indefinite quaternion algebra $B$ over $\Q $.
Assume that the latter holds; we will reach a contradiction.

It follows that $\End _{\Q }(A)\otimes \Q $ is either $\Q $ or a
real or imaginary quadratic field $E$ that embeds in $B$. Indeed,
the possibility $\End _{\Q }(A)\otimes \Q =B$ is not allowed since
otherwise there would exist an embedding $B\hookrightarrow \End
_{\Q }(\Omega ^1_{A})\simeq \M_2(\Q )$, which is not possible
because $B\not \simeq \M_2(\Q )$.

In addition, $\End _{\Q }(A)\otimes \Q$ can not be imaginary
quadratic. Indeed, assume that $\End _{\Q }(A)\otimes \Q =E\subset
B$ were an imaginary quadratic field. Then $B=\End _{\qbar
}(A)=\End _K(A)$ for a certain quadratic number field (cf.
\cite{DiRo}). For any prime $\lambda $ on $E$ over a rational
prime $\ell $, let $\delta _{\lambda }: \mathrm{Gal }(\qbar /\Q
)\rightarrow \Z _{\ell }^*$ be the determinant of the $\lambda
$-adic representation of $\GL _2$-type attached to the $\ell
$-power torsion of $A$ (cf.\,\cite{Ri4}, Section 3) and let $\chi
_{\ell }$ be the $\ell $-adic cyclotomic character. By \cite{Ri4},
Lemma 3.1, there is a character of finite order $\epsilon :
\mathrm{Gal }(\qbar /\Q )\rightarrow E^*$ such that $\delta
_{\lambda } = \epsilon \chi _{\lambda }$, which is trivial if and
only if $K$ is real. By \cite{Ri4}, Lemma 3.2, $\epsilon $ is
even, that is, $\epsilon (\sigma )=+1$ for any complex conjugation
$\sigma $ on $\qbar $. Since $E$ is imaginary, it turns out that
$\epsilon $ is trivial and hence $K$ is real. Then we would obtain
an abelian surface $(A, \iota :B\hookrightarrow \End _K(A)\otimes
\Q )$ that would represent a point on a Shimura curve rational
over $K$. This again yields a contradiction with \cite{Sh}. This
shows that $\End _{\Q }(A)\otimes \Q $ is either $\Q $ or a real
quadratic field. $\Box $

Natural examples of abelian surfaces $A$ with potential
quaternionic multiplication over $\Q $ and $\End _{\Q }(A)$
quadratic arise when considering modular abelian surfaces (with no
character): simple two-dimensional factors of the Jacobian
$J_0(N)$ of the modular curve $X_0(N)$ of level $N\geq 1$. Indeed,
let $A/\Q $ be a modular abelian surface. Although $\End _{\Q
}(A)\otimes \Q =F$ is a quadratic number field, it may be the case
that $\End _{\bar \Q }(A)\otimes \Q $ is a strictly larger algebra
(cf. \cite{Mo}, \cite{Ri4}). This is exactly the case when
$A\stackrel {\Q }{\sim }A_f$ is isogenous to the abelian variety
attached by Shimura to a newform $f\in S_2(\Gamma _0(N))$ without
CM and with an {\em extra-twist}. Then, $\End _{\bar \Q
}(A_f)\otimes \Q $ is necessarily a (possibly splitting)
quaternion algebra $B$ over $\Q $.

More precisely, if we let $\chi $ the quadratic character attached
to the extra-twist of $f$, then $\End _{\qbar }A$ is an order in
the quaternion algebra $(\frac {d,\chi (-1) r}{\Q })$, where $d =
\mbox{ disc }(F)$ and $r$ is the conductor of $\chi $.

Computations due to Hasegawa (\cite{Ha}) and Clark and Stein show
that the only endomorphism algebras $\End _{\bar \Q }(A_f)\otimes
\Q $ that occur for extra-twisting newforms $f$ of level $N<4500$
and trivial Neben-typus are the quaternion algebras $B$ over $\Q $
of discriminant $D = 1$, $6$, $10$, $14$ and $15$.

We note that it very often holds that the ring $\End _{\Q }(A)=
\TT \simeq \Z [\{a_n\}]$ of Hecke operators acting on the modular
abelian surface $A$ is the maximal ring of integers of the number
field $\Q (\{a_n\})$ generated by the Fourier coefficients of $f$.
Under this assumption and the assumptions of Theorem
\ref{ModOrder}, it implies that $\cO =\End _{\qbar }(A)$ contains
a maximal quadratic order and is therefore a {\em primitive}
quaternion order. In particular, $\cO $ is a {\em Bass} order.

\section{Examples of non modular abelian surfaces with potential
quaternionic multiplication}

In this section, we exhibit examples of abelian surfaces with
potential quaternionic multiplication which show that all the
cases of Theorem \ref{ModOrder} can occur.

In order to accomplish that, we consider several particular fibres
of a family of curves of genus two whose Jacobian varieties have
multiplication by a maximal order in the quaternion algebra $B$
over $\Q $ of reduced discriminant $6$ obtained by
Hashimoto-Tsunogai in \cite{HaTs}. We also refer the reader to
\cite{HaMu}. For these curves $C/K$, we compute the minimal field
of definition $L$ of all endomorphisms on $J(C)$ and we apply a
theorem proved by the authors in \cite{DiRo} to conclude that the
varieties $J(C)/K$ are not of $\GL _2$-type over the base field
$K$.

\begin{theorem}

{\bfseries I.} Let $C_1/\Q (\sqrt {2})$ be a smooth projective
model of the curve

$$
Y^2 = (X^2 +5) ( (-1/6 + \sqrt{2}) X^4 + 20 X^3 - 490/6 X^2 + 100
X + 25 (-1/6 - \sqrt{2})).
$$

Then, the Jacobian variety $A_1=J(C_1)/\Q (\sqrt {2})$ of $C_1$
has multiplication by a maximal order $\cO $ in the quaternion
algebra $B_6$ of discriminant $6$ over the quartic extension
$L=\Q(\sqrt{2}, \sqrt{-1}, \sqrt{-5})$ of $\Q (\sqrt {2})$.
Moreover,

\begin{itemize}
\item $\End _{\Q(\sqrt{2}, \sqrt{-5})}(A_1)\otimes \Q =
\Q(\sqrt{2})$

\item $\End _{\Q(\sqrt{2}, \sqrt{5})} (A_1)\otimes \Q =\Q (\sqrt
{3})$

\item $\End _{\Q(\sqrt{2}, \sqrt{-1})}(A_1)\otimes \Q = \Q (\sqrt
{-6})$.

\end{itemize}
and $\End _{\Q (\sqrt {2})}(A_1)=\Z $.

\vspace{0.2cm} {\bfseries II.} Let $C_2/\Q $ be a smooth
projective model of the curve

$$
Y^2 = (X^2 + 7/2) (83/30 X^4 + 14 X^3 - 1519/30 X^2 + 49 X -
1813/120)
$$
and let $A_2=J(C_2)/\Q $ be its Jacobian variety. Then,

\begin{itemize}
\item $\End _{L}(A_2)=\cO $ is a maximal order in $B_6$ for
$L=\Q(\sqrt{-6}, \sqrt{-14})$.

\item $\End _{\Q (\sqrt{-14})}(A_2)\otimes \Q = \Q(\sqrt{2})$

\item $\End _{\Q (\sqrt{21})} (A_2)\otimes \Q =\Q (\sqrt {3})$

\item $\End _{\Q (\sqrt{-6})}(A_2)\otimes \Q = \Q (\sqrt {-6})$

\item $\End _{\Q }(A_2)=\Z $.

\end{itemize}

\end{theorem}

{\em Proof. } By \cite{HaTs}, Lemma 4.5, putting $(\tau _1, \sigma
_1)= (\sqrt{-2}, \sqrt{-1})$ and $(\tau _2, \sigma _2)=
(\sqrt{-3/2}, \sqrt{-3/2})$, we know that $\End _{\qbar }(A_i)
\otimes \Q $, $i=1, 2$ contains the quaternion algebra $B_6$.
Moreover, as it was shown in \cite{HaMu} for an isomorphic family
of curves, $\End _{\qbar }(A_i)$ contains a maximal order in
$B_6$.

The curve $C_1$ is defined over $K=\Q(\sqrt{2})$ and the primes of
bad reduction of the curve all divide $N=30$.

For the primes $\wp $ in $K$ above the rational primes $p=7, 17,
23, 31, 41, 47, 71, 73, 79, \\ 89$ and $97$, we computed the
characteristic polynomial of the image of $\Fr \; \wp$ and we
obtained that the values of the {\em half-traces} $(a_\wp ,
b_\wp)$ are $\pm 2 \sqrt{2}, \pm 2 \sqrt{-6}, \pm 4 \sqrt{2}, \pm
2 \sqrt{3},\\ (8,8), (0,0), \pm 8 \sqrt{3}, \pm 4 \sqrt{-6}, \pm 2
\sqrt{3}, (-16, -16)$ and $(0,0)$, respectively.

Let us remark that the values $a_\wp$ and $b_\wp$ correspond to
the factorization

$$
(x^2 - a_\wp x +p ) (x^2 - b_\wp x +p)
$$
except for the rational primes $p = 17, 73, 97$, namely those that
are a square mod $4$ and a non-square mod $5$. The latter
correspond to the factorization

$$
(x^2 - a_\wp x -p ) (x^2 - b_\wp x -p).
$$

Above, the term $-p$ is due to the fact that $\sqrt{-6}$ is not
real. Indeed, as we explained in \cite{DiRo}, Section 5.1, when an
imaginary quadratic field shows up in the endomorphism algebra,
the Galois representations are reducible but the determinant of
the two-dimensional irreducible components are not $\chi$.

The endomorphism algebra of an abelian surface containing a
quaternion algebra $B$ is either $B$ itself or a matrix algebra
$\M_2(K)$ for $K$ an imaginary quadratic field. Let us explain why
the latter does not occur for $A_1$ and thus $\End_{\qbar } (A_1)$
is a maximal order in $B_6$. To see this, we will show that for a
suitable prime $\ell$, the representation $\rho_\ell$ giving the
action of Galois on the Tate module of $A_1$ is not potentially
abelian, thus eliminating the case of the product of two elliptic
curves with complex multiplication.

The values of the traces computed show that $\End_{K}(A_1)\otimes
\Q $ can not contain a real or imaginary quadratic field. Thus, we
have $\End_{K}(A_1) = \Z$. Hence, the abelian surface $A_1$
acquires multiplication by a quadratic field over a quadratic
extension $K'$ of $K$ unramified outside $30$. Looking at the
values of the half-traces we see that they fall in a fixed
quadratic field only if we restrict to the Galois group of one of
the following quadratic extensions $K'$ of $K$ (among those with
the above specified ramifying places): $K(\sqrt{-1})$,
$K(\sqrt{-5})$ and $K(\sqrt{5})$. Thus we know that over one of
these fields, $\End_{K'}(A_1)\otimes \Q $ contains a quadratic
field and in particular the representations $\rho_\ell$ restricted
to the corresponding Galois group decompose and have two
two-dimensional irreducible components:
$$
\rho_\ell|_H \cong \sigma_\ell \oplus \sigma'_\ell
$$
where $H= \Gal(\Q / K')$ and $K'$ is one of above the three
fields.

The next and last step is to check that for a suitable $\ell$ the
images of these  two-dimensional components are maximal. We have
checked this for the three possibilities for $K'$.  Here we give
the details only for the case $K' = K(\sqrt{5})$, the other two
being similar. In this case, when we restrict to $H$ the traces of
the two-dimensional components fall in $\Q(\sqrt{3})$, so we are
assuming that $A_1$ has real multiplication by this field, defined
over $K(\sqrt{5})$. We take the prime $\ell = 11$ and observe that
it decomposes in $\Q(\sqrt{3})$. Let us show that the image of
$\sigma_{11}$ is the full $\GL(2, \Z_{11})$. Thanks to a
well-known lemma of Serre, it is enough to show that the residual
representation $\bar{\sigma}_{11}$ has maximal image $\GL(2,
\F_{11})$. There is a well-known procedure to prove this
maximality computationally, we recall it for the reader's
convenience: first, compute the reduction modulo $11$ of the
characteristic polynomials of several Frobenius elements
(corresponding to elements in $H$), and observe that some have two
different roots in $\F_{11}$, one of them of maximal order $\ell-1
= 10$, and on the other hand some have irreducible characteristic
polynomials whose roots are elements of high order in
$\F_{121}^*$. We have checked this using the Frobenius elements
for the following primes (all splitting totally in $K(\sqrt{5})$):
$p = 31, 41, 71, 79$ and $89$. Using the classification of maximal
subgroups of $\GL_2(\F_\ell)$  due to L. E. Dickson, one concludes
that the image of $\bar{\sigma}_{11}$ is the full $\GL(2,
\F_{11})$. Therefore, the image of $\sigma_{11}$ is also maximal.
The same holds if we take the other two quadratic extensions $K'$
of $K$ mentioned above. Thus, we conclude that $A_1$ is not in the
potentially CM case, and so it has only potential quaternionic
multiplication.

We know that $\End_{\qbar}(A_1)$ is a maximal order in $B_6$, and
that $\End_K(A_1) = \Z$. In consequence, Theorem 1.3 in [DiRo]
yields that $[L:K]=4$. By considering all possible quadratic
extensions unramified outside $\{ 2,3,5 \}$ and by applying Lemma
5.2 in \cite{DiRo}, we conclude that the traces above only match
with the fields of definition and intermediate endomorphism
algebras claimed in part {\bfseries I} of our theorem. A similar
argument and similar computations yield our statement for curve
$C_2$. $\Box $

\begin{corollary}

The Jacobian variety $J(C_2)/\Q $ is a non modular abelian surface
with potential quaternionic multiplication over $\Q $.

\end{corollary}

Note that, since $B_6=(\frac {-6, 2}{\Q })=(\frac {-6, 3}{\Q })$,
the three possible intermediate endomorphism algebras $\Q (\sqrt
{2})$, $\Q (\sqrt {3})$, $\Q (\sqrt {-6})\subset B_6$ allowed by
the Theorem 1.3 in \cite{DiRo} arise both on $A_1$ and $A_2$.

\section{Fake elliptic curves over quadratic imaginary fields}

Finally, let us conclude by exhibiting a pair of examples of
abelian surfaces $A/K$ over a quadratic imaginary field $K$ such
that all quaternionic endomorphisms of $A$ are defined over $K$
itself. These are therefore examples of what we can call fake
elliptic curves over $K$. Our examples are non trivial, since it
can be shown that they are not the base extension to $K$ of an
abelian surface over $\Q $.

As a by-product, we also show that all computations performed in
\cite{HaTs} supporting an analogue of the Sato-Tate conjecture for
these surfaces are unconditionally correct (cf. \cite{HaTs} for
details).

In order to accomplish that, assume that for a particular curve
$C/K$, computations suggest that the minimal field of definition
of the ring of correspondences of $C$ is $L=K$. This can be the
case if the characteristic polynomial of $\Fr \; \wp$, for $\wp$ a
prime in $K$ of good reduction of $A$ and residual degree $1$,
factorizes as

$$
(x^2 - a_\wp x + p)^2 \quad \quad \quad \quad \eqno(2.1)
$$
with $a_\wp \in \Z$.

In order to eliminate the quadratic case $\Gal(L/K) = C_2$ allowed
by \cite{DiRo}, Theorem 1.3, we may compute all quadratic
extensions $L$ of $K$ ramifying only at the primes of bad
reduction of $A$ and, for each of them, exhibit a prime $\wp$ of
$K$ inert in this quadratic extension verifying $(2.1)$ with
$a_\wp \neq 0$. This contradicts formula (4.5) in \cite{DiRo} and
we thus conclude that the field $L$ of definition of the
quaternionic multiplication can not be a quadratic extension of
$K$.

Moreover, applying \cite{DiRo}, Lemma 5.2, we see that the above
fact is also incompatible with the case $\Gal(L/K)=D_2$ because it
violates the trace $0$ condition for those primes that do not
totally decompose in $L/K$. As examples of this phenomenon we can
exhibit the following: in \cite{HaTs}, (3.1), a family $S_6(t,s)$
of QM-curves $C_{(t,s)}$ of genus $2$ is given. It is such that,
for every rational value of the parameter $t$, the curve
$C_{(t,s)}$ is defined over the imaginary quadratic field
$\Q(s)=\Q(\sqrt{-3 + 14 t^2 - 27 t^4})$. For several rational
values of $t$, the authors considered the action of
$\Gal(\qbar/K)$ on the Tate modules of $J(C_{(t,s)})$ and they
observed that, for the first $30$ primes of residual degree $1$ in
$K$, formula (2.1) is verified. In consequence, they suggested
that all endomorphisms are defined over $K$, that is, $L=K$. In
fact, the experimental verification of the Sato-Tate conjecture
that they obtain in their article depends on this assumption (cf.
the examples 1 and 2 in \cite{HaTs}, page 1655, corresponding to
$(t,s)= (2, \sqrt{-379})$ and $(2/3, \sqrt{-19}/3)$, respectively,
and the tables and figures of p. 1658, 1659).

For these two particular examples, we computed all quadratic
extensions of $K$ ramifying only at primes of bad reduction and
for each of them we found a prime $\wp$ of $K$ inert in this
extension verifying (2.1) with $a_\wp \neq 0$. This eliminates all
the cases listed in \cite{DiRo}, Theorem 1.3, except $L=K$. We
thus conclude

\begin{theorem}

Let $C_3/K_3$, $C_4/K_4$ be the fibres of the Hashimoto-Tsunogai's
family $S_6(t,s)$ at the values $t_3= 2$ and $t_4=2/3$ over $K_3 =
\Q(\sqrt{-379})$ and $K_4=\Q(\sqrt{-19})$, respectively. Then,
$\End _{K_i}(J(C_i))\simeq \cO $ is a maximal order in $B_6$.

\end{theorem}

In particular, all computations performed in \cite{HaTs} for these
two examples supporting the Sato-Tate conjecture for these
surfaces are unconditionally correct.

\begin{remark} The computation of the absolute Igusa invariants of
these curves show that there exists no curve $C/\Q $ such that
$C\simeq C_3$ nor $C\simeq C_4$ over $\bar \Q $. See \cite{Me} for
details.
\end{remark}

\begin{remark} The importance of the property $L=K$ is that it
implies that (2.1) holds for every prime and it is then enough to
compute the number of rational points on the curve over $\F _p$ in
order to determine the characteristic polynomial of $\Fr \; \wp$
for any prime $\wp \mid p $ in $K$, while in general it is
necessary to compute also the number of rational points over
$\F_{p^2}$. This enormously speeds the computations.
\end{remark}

\begin{remark} The above two curves $C_3$ and $C_4$ and also the
respective rings of correspondences on them are defined over
completely imaginary fields $L=K$. This is not a coincidence since
there do not exist curves $C/K$ of genus $2$ over a number field
$K$ admitting a real archimedean place such that $\End _{K}(J(C))$
is a quaternion order. Indeed, this follows from Shimura's result
that the set of real points on Shimura curves is the empty set
(\cite{Sh}).
\end{remark}

\end{document}